\documentclass[11pt]{amsart}
\usepackage{amsmath,amssymb,amsthm,booktabs}
\usepackage[margin=1in]{geometry}
\newcommand{\Rsub}[1]{\ifnum 0<0#1\relax\the\numexpr#1-1\relax\else#1-1\fi}
\newcommand{\R}[2]{R^{#1}_{\Rsub{#1}}(C_{#2})}
\newcommand{\RP}[2]{R^{#1}_{\Rsub{#1}}(P_{#2})}     
\newcommand{\circum}{\operatorname{circ}}      
\newcommand{\fl}[1]{\left\lfloor #1 \right\rfloor}
\newcommand{\cl}[1]{\left\lceil #1 \right\rceil}

\theoremstyle{plain}
\newtheorem{theorem}{Theorem}[section]
\newtheorem{lemma}[theorem]{Lemma}
\newtheorem{proposition}[theorem]{Proposition}
\newtheorem{corollary}[theorem]{Corollary}

\theoremstyle{definition}
\newtheorem{definition}[theorem]{Definition}

\theoremstyle{remark}

\title{The $(t-1)$-chromatic Ramsey number of cycles for three colours}
\author{A. Khamseh and M. H. Safareh}
\address{amir khamseh \\ department of mathematics, kharazmi university}
\email{khamseh@khu.ac.ir; amirkhamseh114@gmail.com}
\address{mohammad hasan safareh \\ department of mathematics, kharazmi university}
\email{mh.safareh@khu.ac.ir}

\date{}

\begin{document}

\begin{abstract}
For integers $t\ge2$ and $\ell\ge3$, let $R^t_{t-1}(C_\ell)$ be the least integer $n$ such that
every colouring of the edges of $K_n$ with $t$ colours contains a cycle of length $\ell$ whose edges
use at most $t-1$ colours. Buci\'{c} and Khamseh determined the analogous number for paths and
asked for its value for cycles. We answer this question for three colours: we show that
$R^3_2(C_\ell)=\ell+\lfloor(\ell+1)/6\rfloor$ for every $\ell\ge11$, and we determine the values
for $3\le\ell\le10$. 
\end{abstract}

\maketitle

\section{Introduction}\label{sec:intro}

Ramsey theory rests on the principle, made precise by Ramsey's theorem~\cite{Ram29}, that every
sufficiently large structure contains a highly ordered substructure. For graphs $G_1,\dots,G_t$ the
classical \emph{$t$-colour Ramsey number} is the least $n$ such that every colouring of the edges of
$K_n$ with $t$ colours contains, for some $i$, a copy of $G_i$ in colour $i$. Exact values are known
only for a handful of families. Among the earliest are those for paths, due to Gerencs\'er and
Gy\'arf\'as~\cite{GG67}, for stripes (disjoint copies of $P_2$), due to Cockayne and
Lorimer~\cite{CL75}, and for arbitrary linear forests, due to Faudree and Schelp~\cite{FS76}. For
three or more colours the picture is far less complete. {\L}uczak~\cite{Luc99} proved
$R(C_n,C_n,C_n)\le(4+o(1))n$, and the three-colour Ramsey number of paths was determined exactly for
long paths by Gy\'arf\'as, Ruszink\'o, S\'ark\"ozy and Szemer\'edi~\cite{GRSS07}; for four or more
colours the problem remains open, the best known lower and upper bounds for paths and even cycles
being due to Yongqi, Yuansheng, Xu and Bingxi~\cite{YYXB06} and to Knierim and Su~\cite{KS19},
respectively. The diagonal multicolour problem for complete graphs has also seen dramatic recent
progress: Balister, Bollob\'as, Campos, Griffiths, Hurley, Morris, Sahasrabudhe and
Tiba~\cite{BBCGHMST26} obtained, for every fixed number of colours $r\ge3$, the first exponential
improvement on the Erd\H{o}s--Szekeres upper bound. Many variants of Ramsey numbers of paths and
cycles have been studied, among them list Ramsey numbers~\cite{ABKKS21}, bipartite
versions~\cite{BLS19a,BLS19d,DGKRS20,WZHT21}, directed and oriented
versions~\cite{BLS19b,BLS19c,LS20}, star-critical and critical-graph
refinements~\cite{HI11,Hoo15}, Ramsey goodness~\cite{AMM25}, Ramsey indexes~\cite{CZ24} and Ramsey
numbers of trees~\cite{MPY25}; we refer to the survey of Conlon, Fox and Sudakov~\cite{CFS15}.

The generalisation studied here was introduced by Erd\H{o}s, Hajnal and Rado~\cite{EHR65}. For a graph
$G$ and integers $t\ge k\ge1$, let $R^t_k(G)$ be the least $N$ such that every $t$-colouring of the
edges of $K_N$ contains a copy of $G$ whose edges use at most $k$ of the $t$ colours. The case $k=1$ is
the classical Ramsey number and is the hardest variant, while $k=t-1$ is the first non-trivial case
and has attracted considerable independent interest. Erd\H{o}s, Hajnal and Rado treated only
$G=K_n$; Erd\H{o}s and Szemer\'edi~\cite{ES72} showed that $R^t_{t-1}(K_n)$ grows exponentially in
$n$, the main open problem being the dependence of the exponent on $t$, and the best known general
bounds are due to Alon, Erd\H{o}s, Gunderson and Molloy~\cite{AEGM02}. In a closely related spirit,
Alon, Buci\'{c}, Christoph and Krivelevich~\cite{ABCK24} determined, up to a logarithmic factor,
how large a connected (and more generally $k$-connected) subgraph using at most $s$ colours can be
guaranteed in any $r$-colouring of $K_n$, which resolves a generalisation, due to Bollob\'as, of the
classical connected-component problem of Gerencs\'er and Gy\'arf\'as~\cite{GG67}. For a comprehensive
account of such ``weakened'' Ramsey-type parameters, including weakened Ramsey, Gallai--Ramsey, Schur
and Rado numbers, we refer to the monograph of Budden~\cite{Bud26}; recent progress on the
star-critical refinement is in Budden, Moun and Jakhar~\cite{BMJ25}.

Let $t\ge 2$. A \emph{$t$-colouring} of $K_n$ is an arbitrary map $\phi\colon E(K_n)\to[t]=\{1,\dots,t\}$.
A subgraph of $K_n$ is \emph{colour-avoiding} if at least one of the $t$ colours does not appear on
its edges. For a graph $H$, the \emph{$(t-1)$-chromatic Ramsey number} $R^t_{t-1}(H)$ is the least
$n$ such that every $t$-colouring of $K_n$ contains a colour-avoiding copy of $H$. Chung and
Liu~\cite{CL78} (see also~\cite{CCL77}) were the first to study this number for non-complete graphs,
and it has since been investigated for stars, matchings, complete graphs and
designs~\cite{BB18,BSZ17,Gya18,GS20,HM99,Jac82,KO17,KO22}, and, closest to the present work, for
paths and linear forests~\cite{BK23,Kha18,Kha20,Kha25,KO12}. For the path $P_\ell$ on $\ell$
vertices, Buci\'{c} and Khamseh~\cite[Theorem~1.1]{BK23} proved that
\begin{equation}\label{eq:BK}
  \RP{t}{\ell}=\ell+\fl{\frac{\ell-2}{2^t-2}}\qquad(t\ge 2,\ \ell\ge 2),
\end{equation}
and asked the following question~\cite[Section~4]{BK23}, which we quote verbatim.

\begin{quote}
``\textbf{Question 4.2.} What is the smallest $n$ such that in any $t$-edge colouring of $K_n$ we can
find a cycle of length $\ell$ using at most $t-1$ colours?''
\end{quote}

That is, Question~4.2 asks for $\R{t}{\ell}$, where $C_\ell$ is the cycle of length $\ell$. We
answer it for three colours.

\begin{theorem}\label{thm:main}
For every integer $\ell\ge3$,
\[
  \R{3}{\ell}=
  \begin{cases}
    \ell+\fl{\dfrac{\ell+1}{6}} & \text{if } \ell\ge 11,\\[2ex]
    \ell+2 & \text{if } \ell=3 \text{ or } 9\le\ell\le 10,\\[1ex]
    \ell+1 & \text{if } 4\le\ell\le 8.
  \end{cases}
\]
\end{theorem}

\begin{proof}
For $3\le\ell\le 11$ this is Proposition~\ref{prop:basic-t3}, and for $\ell\ge11$ it is
Theorem~\ref{thm:t3all}; for $\ell=11$ both give $13$.
\end{proof}

Since $C_\ell$ contains $P_\ell$, \eqref{eq:BK} gives $\R{3}{\ell}\ge\RP{3}{\ell}=\ell+\fl{(\ell-2)/6}$; the first length at which
$\R{3}{\ell}$ exceeds both $\ell+2$ and the path value is $\ell=17$, where $\R{3}{17}=20$ and
$\RP{3}{17}=19$. Section~\ref{sec:basic} reduces the problem to the circumferences of the three
graphs obtained by deleting one colour and settles $\ell\le11$; Section~\ref{sec:t3all} treats
$\ell\ge11$. For $t\ge4$ colours Question~4.2 remains open.

\section{Preliminaries}\label{sec:prelim}

\subsection{Notation}
Graphs are finite and simple. For a graph $G$ we write $|G|$ for the number of vertices, $e(G)$ for
the number of edges, $N_G(v)$ and $\deg_G(v)=|N_G(v)|$ for the neighbourhood and degree of $v$,
$G-v$ for the graph obtained by deleting $v$, and $G+uv$ for the graph obtained by adding the edge
$uv$. For disjoint vertex sets $A,B$, $e_G(A,B)$ is the number of edges of $G$ between $A$ and $B$.
For integers $a\le b$, $[a,b]:=\{a,\dots,b\}$. A path has \emph{order} equal to its number of
vertices; $P_k$ is the path of order $k$ and $C_k$ the cycle of length $k$. The
\emph{circumference} $\circum(G)$ is the largest length of a cycle in $G$ ($0$ if $G$ is acyclic).
A Hamilton cycle (path) of $G$ is a cycle (path) through all vertices of $G$.

\emph{Blocks.} A \emph{block} of $G$ is a maximal connected subgraph without a cut vertex (of
itself); we only consider blocks with at least two vertices, that is, maximal $2$-connected
subgraphs and bridges with their ends. We use the following standard facts~\cite[Section~3.1]{Die17}.
\begin{itemize}
\item[(B1)] Two distinct blocks share at most one vertex, and every edge lies in exactly one block.
\item[(B2)] Every cycle lies in a block (a cycle is $2$-connected); also $\circum(B)\le\circum(G)$ for
  every block $B$, as $B\subseteq G$.
\item[(B3)] Every block $B$ is an induced subgraph: if $x,y\in V(B)$ are adjacent in $G$, then $B+xy$
  is connected without a cut vertex, so $xy\in E(B)$ by maximality.
\item[(B4)] $\sum_B(|B|-1)=|G|-c(G)$, where $c(G)$ is the number of components and the sum is over all
  blocks. For connected $G$ this follows by induction on the number of blocks from
  \cite[Lemma~3.1.4]{Die17} (the block graph of a connected graph is a tree), deleting from $G$ all
  vertices of a leaf block except its cut vertex; in general, sum over the components.
\end{itemize}

\emph{Colourings.} Let $\phi$ be a $t$-colouring of $K_n$. For $j\in[t]$ the \emph{$j$-avoiding
graph} $G_j=G_j(\phi)$ is the spanning subgraph of $K_n$ with edge set $\{e:\phi(e)\ne j\}$. A
subgraph of $K_n$ is colour-avoiding if and only if it lies in some $G_j$. Every edge lies in exactly
$t-1$ of the graphs $G_j$, so
\begin{equation}\label{eq:prelim-sum}
  \textstyle\sum_{j\in[t]}e(G_j)=(t-1)\binom n2 .
\end{equation}
Restricting a $t$-colouring of $K_{n'}$ to $n\le n'$ vertices shows that if every $t$-colouring of
$K_n$ has a colour-avoiding $C_\ell$, then so does every $t$-colouring of $K_{n'}$. Hence
$\R{t}{\ell}\le n$ if and only if every $t$-colouring of $K_n$ has a colour-avoiding $C_\ell$
(finiteness of $\R{3}{\ell}$ follows from the upper bounds proved below), and trivially
$\R{t}{\ell}\ge\ell$.

\subsection{Cited results}

\begin{theorem}[Dirac~\cite{Dir52}]\label{thm:dirac}
Every $2$-connected graph $G$ with minimum degree at least $k$ satisfies $\circum(G)\ge\min\{|G|,2k\}$.
\end{theorem}

\begin{theorem}[Bondy--Chv\'atal~\cite{BC76}]\label{thm:closure}
Let $H$ be a graph on $N$ vertices. Repeatedly joining two non-adjacent vertices whose degree sum in
the current graph is at least $N$, until no such pair remains, yields a graph $\mathrm{cl}(H)$, the
\emph{closure} of $H$, which does not depend on the order in which the edges are added.
\end{theorem}

\begin{theorem}[Woodall~\cite{Woo72}; see~{\cite[Theorem~8]{LN23}}]\label{thm:woodall}
Let $k\ge 0$ and let $G$ be a graph on $n\ge 2k+3$ vertices with
$e(G)\ge\binom{n-k-1}{2}+\binom{k+2}{2}+1$. Then $G$ contains a cycle of length $m$ for every
$m\in[3,n-k]$.
\end{theorem}

\begin{theorem}[Kopylov~\cite{Kop77}; see~{\cite[Theorem~1.3]{FKLV}}]\label{thm:kopylov}
Let $n\ge k\ge 5$ and let $G$ be a $2$-connected graph on $n$ vertices with $\circum(G)<k$. Then
$e(G)\le\max\{h(n,k,2),h(n,k,\fl{(k-1)/2})\}$, where $h(n,k,a):=\binom{k-a}2+a(n-k+a)$.
\end{theorem}

\begin{theorem}[Gerencs\'er--Gy\'arf\'as~\cite{GG67}]\label{thm:gg}
Let $a\ge b\ge 2$ and $N\ge a+\fl{b/2}-1$. Every red/blue colouring of $K_N$ contains a red path of
order $a$ or a blue path of order $b$.
\end{theorem}

\section{lower bounds and small lengths}\label{sec:basic}

From now on $t=3$. By~\eqref{eq:prelim-sum}, every $3$-colouring of $K_n$ satisfies
\begin{equation}\label{eq:basic-sum}
  e(G_1)+e(G_2)+e(G_3)=n(n-1).
\end{equation}


\begin{lemma}[]\label{lem:short}
Let $c\ge 5$ and let $\phi$ be a $3$-colouring of $K_n$. If some $G_j$ contains a cycle $C$ of
length $c$, then some $G_i$ contains a cycle of length $c-1$ with vertex set contained in $V(C)$.
\end{lemma}

\begin{proof}
Let $C=v_0v_1\cdots v_{c-1}v_0$, indices modulo $c$. A cycle \emph{avoids} $i$ if it lies in $G_i$.
Suppose that no cycle of length $c-1$ on vertices of $C$ avoids any colour.

\emph{Step 1.} Every chord $v_av_{a+2}$ (not an edge of $C$, as $c\ge5$) has colour $j$: otherwise
$v_av_{a+2}v_{a+3}\cdots v_{a-1}v_a$ is a cycle of length $c-1$ avoiding $j$.

\emph{Step 2: $c$ odd.} As $2$ is invertible modulo $c$, $W=v_0v_2v_4\cdots v_{2(c-1)}v_0$ is a cycle
through $V(C)$; its edges are the pairs $v_xv_y$ with $y-x\equiv\pm2$, all of colour $j$ by Step~1.
The pair $g=v_0v_4$ is not an edge of $W$, since $4\not\equiv\pm2\pmod c$ for odd $c\ge5$. Take
$i\notin\{j,\phi(g)\}$. Deleting $v_2$ from $W$ and adding $g$ gives a cycle of length $c-1$ avoiding
$i$, a contradiction.

\emph{Step 3: $c=2m$ even, so $m\ge3$.} The cycles $D_0=v_0v_2\cdots v_{2m-2}v_0$ and
$D_1=v_1v_3\cdots v_{2m-1}v_1$ have all their edges, the \emph{$D$-edges}, of colour $j$. Call
$v_xv_y$ with $x\not\equiv y\pmod2$ a \emph{cross pair}. For $a,b$ of different parity, let
$Z(a,b)$ run along $D_{a\bmod2}$ from $v_{a+2}$ forwards to $v_a$, then along $D_{b\bmod2}$ from $v_b$
backwards to $v_{b+4}$, and close up; and let $Z'(a,b)$ run along $D_{a\bmod2}$ from $v_{a+2}$
forwards to $v_a$, then from $v_{b+4}$ forwards to $v_b$, and close up. Both are cycles on
$m+(m-1)=c-1$ vertices (they omit $v_{b+2}$), consisting of $D$-edges and two cross pairs:
$v_av_b,v_{b+4}v_{a+2}$ for $Z(a,b)$ and $v_av_{b+4},v_bv_{a+2}$ for $Z'(a,b)$. Such a cycle avoids
a colour $i\ne j$ unless its two cross pairs have the two colours of $[3]\setminus\{j\}=:\{\alpha,
\beta\}$. Hence every cross pair $v_av_b$ has colour in $\{\alpha,\beta\}$, and the symmetric map
$\psi(x,y):=\phi(v_xv_y)\in\{\alpha,\beta\}$ ($x,y$ of different parity) satisfies
\begin{equation}\label{eq:basic-R}
  \psi(a,b)\ne\psi(a+2,b+4)\quad\text{(from $Z(a,b)$)},\qquad
  \psi(a,b+4)\ne\psi(a+2,b)\quad\text{(from $Z'(a,b)$)}.
\end{equation}
As $\psi$ is two-valued, applying the first relation to $(x+2,y-4)$ and the second to $(x,y-4)$ gives
$\psi(x+4,y)=\psi(x,y)$, and by symmetry of $\psi$ also $\psi(x,y+4)=\psi(x,y)$. With the first
relation this yields $\psi(a+2,b)=\psi(a+2,b+4)\ne\psi(a,b)$, and by symmetry
\begin{equation}\label{eq:basic-alt}
  \psi(a+2,b)\ne\psi(a,b)\qquad\text{and}\qquad\psi(a,b+2)\ne\psi(a,b).
\end{equation}
Iterating~\eqref{eq:basic-alt} $m$ times and using $v_{a+2m}=v_a$ shows that $m$ is even, a
contradiction if $m$ is odd. Let $m$ be even, so $m\ge4$. The pair $g=v_0v_4$ is not a $D$-edge, as
$4\not\equiv\pm2\pmod{2m}$. Choose $i\in\{\alpha,\beta\}\setminus\{\phi(g)\}$ and let $\pi$ be the
other element of $\{\alpha,\beta\}$. By~\eqref{eq:basic-alt}, $\psi(4,1)\ne\psi(4,3)$, so
$\psi(4,y)=\pi$ for some odd $y$, and $\psi(6,y-2)\ne\psi(4,y-2)\ne\psi(4,y)$ gives
$\psi(6,y-2)=\pi$. The cycle
\[
  v_6v_8\cdots v_{2m-2}v_0\,v_4\,v_yv_{y+2}\cdots v_{y-2}\,v_6
\]
(all of $D_0$ except $v_2,v_4$ in order along $D_0$, then $v_4$, then all of $D_1$ along $D_1$) has
$(m-1)+m=c-1$ vertices, and its edges are $D$-edges, $g$ and the cross pairs $v_4v_y$, $v_{y-2}v_6$
of colour $\pi$. So it avoids $i$, a contradiction.
\end{proof}

\begin{corollary}\label{cor:mono}
Let $\ell\ge 4$. A $3$-colouring of $K_n$ has a colour-avoiding $C_\ell$ if and only if
$\max_j\circum(G_j)\ge\ell$. Consequently $\R{3}{\ell}\le\R{3}{\ell+1}$.
\end{corollary}

\begin{proof}
If $\max_j\circum(G_j)\ge\ell$, choose $j$ and a cycle $Z$ of $G_j$ with $|Z|\ge\ell$ and $|Z|$
minimal. If $|Z|>\ell$, then $|Z|\ge5$ and Lemma~\ref{lem:short} gives a cycle of length
$|Z|-1\ge\ell$ in some $G_i$, contradicting minimality. The converse is trivial. If every
$3$-colouring of $K_n$ has a colour-avoiding $C_{\ell+1}$, it has $\max_j\circum(G_j)>\ell$, hence
a colour-avoiding $C_\ell$.
\end{proof}

\subsection{Lower bounds}

\begin{lemma}[]\label{lem:basic-gap}
Let $V(H)=S\cup T\cup I$ be a partition with $S\ne\emptyset$, $I$ independent and no edges between $I$
and $T$. Then $\circum(H)\le2|S|$ if $T=\emptyset$, and $\circum(H)\le2|S|+|T|-1$ otherwise.
\end{lemma}

\begin{proof}
Let $Z$ be a cycle and $s'=|V(Z)\cap S|$. Neighbours of $I$-vertices lie in $S$. If $s'=0$, then
$V(Z)\subseteq T$, so $T\ne\emptyset$ and $|Z|\le|T|\le2|S|+|T|-1$. If $s'\ge1$, the $S$-vertices cut
$Z$ into $s'$ \emph{gaps} (possibly empty segments between cyclically consecutive $S$-vertices). An
$I$-vertex on $Z$ forms a gap by itself, so the gaps containing $T$-vertices contain no $I$-vertex.
Hence $|V(Z)\cap I|\le s'-1$ if $Z$ meets $T$, and $|V(Z)\cap I|\le s'$ otherwise. In the first case
$|Z|\le 2s'-1+|T|$, in the second $|Z|\le 2s'$; both are at most the claimed bound.
\end{proof}

\begin{proposition}\label{prop:lowA}
For every $\ell\ge3$, $\R{3}{\ell}\ge\ell+\fl{(\ell+1)/6}$.
\end{proposition}

\begin{proof}
For $\ell\le4$ the bound is $\R{3}{\ell}\ge\ell$. Let $\ell\ge5$, $d:=\fl{(\ell+1)/6}\ge1$ and
$n:=\ell-1+d$. Partition $V(K_n)=A_1\cup A_2\cup A_3$ with $|A_1|=\ell-2d$, $|A_2|=2d-1$, $|A_3|=d$
(all positive, as $2d\le(\ell+1)/3<\ell$). Colour edges inside $A_j$ with $j$, and edges between $A_i$
and $A_k$, $i<k$, with $k$. In $G_3$ the vertices of $A_3$ are isolated, so
$\circum(G_3)\le|A_1\cup A_2|=\ell-1$. In $G_1$ the set $A_1$ is independent, and
Lemma~\ref{lem:basic-gap} with $(I,S,T)=(A_1,A_2\cup A_3,\emptyset)$ gives $\circum(G_1)\le
2(3d-1)\le\ell-1$. In $G_2$ the set $A_2$ is independent with no edges to $A_1$, and
Lemma~\ref{lem:basic-gap} with $(I,S,T)=(A_2,A_3,A_1)$ gives $\circum(G_2)\le2d+\ell-2d-1=\ell-1$. So
there is no colour-avoiding $C_\ell$, and $\R{3}{\ell}>n$.
\end{proof}

\begin{proposition}\label{prop:lowC}
For every $\ell\ge9$, $\R{3}{\ell}\ge\ell+2$.
\end{proposition}

\begin{proof}
Let $n=\ell+1\ge10$ and $V(K_n)=\{b_1,b_2,c_1,c_2\}\cup A$. Give colour $2$ to $b_1b_2$, $b_1c_1$,
$b_2c_2$ and to all edges between $\{b_1,b_2\}$ and $A$; colour $3$ to $c_1c_2$, $b_1c_2$, $b_2c_1$ and
to all edges between $\{c_1,c_2\}$ and $A$; colour $1$ to the edges inside $A$. Then $b_1,b_2$ have
degree $1$ in $G_2$ (their only edges not of colour $2$ are $b_1c_2$, $b_2c_1$) and $c_1,c_2$ have
degree $1$ in $G_3$ (only $c_1b_1$, $c_2b_2$), so $\circum(G_2),\circum(G_3)\le n-2$. In $G_1$ the set
$A$ is independent, and Lemma~\ref{lem:basic-gap} with $(I,S,T)=(A,\{b_1,b_2,c_1,c_2\},\emptyset)$ gives $\circum(G_1)\le8\le n-2$.
So no $G_j$ contains $C_{n-1}=C_\ell$, and $\R{3}{\ell}>n=\ell+1$.
\end{proof}

\begin{lemma}\label{lem:basic-K4}
$\R{3}{3}\ge 5$ and $\R{3}{4}\ge 5$.
\end{lemma}

\begin{proof}
Colouring the three perfect matchings of $K_4$ with distinct colours gives every triangle three
colours, as its edges are pairwise adjacent. Colouring $ab,ac$ with $1$, $bc,bd$ with $2$ and $cd,da$ with $3$ makes every colour class two
adjacent edges; a $4$-cycle avoiding $j$ would leave both edges of colour $j$ in the complementary
perfect matching, which is impossible.
\end{proof}

\subsection{Small lengths}

\begin{proposition}\label{prop:woodall}
Let $\ell\ge 3$ and $\ell\le n\le 2\ell-3$ with
\begin{equation}\label{eq:basic-wood}
  n(n-1)>3\Bigl[\tbinom{\ell-1}{2}+\tbinom{n-\ell+2}{2}\Bigr].
\end{equation}
Then $\R{3}{\ell}\le n$.
\end{proposition}

\begin{proof}
By~\eqref{eq:basic-sum} and integrality, some $G_j$ has at least
$\binom{\ell-1}{2}+\binom{n-\ell+2}{2}+1$ edges. Theorem~\ref{thm:woodall} with $k=n-\ell\ge0$
applies, as $n\ge2k+3$ is $n\le2\ell-3$, and gives $C_\ell\subseteq G_j$ since $\ell\in[3,n-k]$.
\end{proof}

For $(\ell,n)=(8,9)$, \eqref{eq:basic-wood} fails by equality ($72=3(21+3)$); the next lemma
handles this case.

\begin{lemma}\label{lem:basic-nine}
Let $G$ be a graph on $9$ vertices with $e(G)\ge 24$ and $\circum(G)\le 7$. Then $G$ consists of a
$K_7$ and a $K_3$ sharing one vertex. In particular one vertex has degree $8$ and all others have
degree $6$ or $2$.
\end{lemma}

\begin{proof}
If $G$ were $2$-connected, Theorem~\ref{thm:kopylov} with $n=9$ and $k=8$ would give
$e(G)\le\max\{h(9,8,2),\,h(9,8,3)\}=\max\{21,\,22\}<24$. Let $B_1,\dots,B_r$ be the blocks and $x_i:=|B_i|-1\ge1$;
by (B1) and (B4), $e(G)=\sum_ie(B_i)$ and $\sum_ix_i\le8$, and by (B2) $\circum(B_i)\le7$. No block
has $9$ vertices (by (B3) it would equal $G$, which is not $2$-connected). If $|B_1|=8$, then Theorem~\ref{thm:kopylov} with
$n=k=8$ gives $e(B_1)\le\max\{h(8,8,2),h(8,8,3)\}=19$, the other blocks have $\sum_{i\ge2}x_i\le1$,
and $e(G)\le20$. So $x_i\le6$ for all $i$, and $e(B_i)\le x_i(x_i+1)/2$. If $x_1=6$, then
$\sum_{i\ge2}x_i\le2$ and the other blocks have at most $3$ edges, with equality only for a single
$K_3$; so $e(G)\le 24$ with equality only if $G$ has exactly the two blocks $K_7$ and $K_3$. If all
$x_i\le5$, then $x_i(x_i+1)/2\le3x_i$ with equality only for $x_i=5$, so $e(G)\le24$ with equality
only if $\sum x_i=8$ is a sum of $5$'s, which is impossible. Hence $G$ has the two blocks $K_7$ and
$K_3$ and $\sum x_i=8$, so $G$ is connected by (B4), every vertex lies in a block, and as $7+3=10$
the blocks share exactly one vertex by (B1). The degrees follow.
\end{proof}

\begin{proposition}\label{prop:basic-t3}
For $3\le\ell\le 11$ the values of $\R{3}{\ell}$ are as in Table~\ref{tab:basic-t3}.
\end{proposition}

\begin{table}[ht]
\centering
\begin{tabular}{@{}lccccccccc@{}}
\toprule
$\ell$ & 3 & 4 & 5 & 6 & 7 & 8 & 9 & 10 & 11\\
$\R{3}{\ell}$ & 5 & 5 & 6 & 7 & 8 & 9 & 11 & 12 & 13\\
lower bound & \ref{lem:basic-K4} & \ref{lem:basic-K4} & \ref{prop:lowA} & \ref{prop:lowA} &
  \ref{prop:lowA} & \ref{prop:lowA} & \ref{prop:lowC} & \ref{prop:lowC} & \ref{prop:lowC}\\
upper bound & matchings & \ref{prop:woodall} & \ref{prop:woodall} & \ref{prop:woodall} &
  \ref{prop:woodall} & \ref{lem:basic-nine} & \ref{prop:woodall} & \ref{prop:woodall} &
  \ref{prop:woodall}\\
\bottomrule
\end{tabular}
\caption{Values of $\R{3}{\ell}$ for $3\le\ell\le 11$ and the results giving the lower and upper
bounds.}
\label{tab:basic-t3}
\end{table}

\begin{proof}
The lower bounds are given in the table. 

\emph{$\ell=3$:} if a $3$-colouring of $K_5$ has no
colour-avoiding triangle, two edges $xy,xz$ of the same colour cannot exist (the triangle $xyz$
would miss a colour), so each colour class is a matching with at most $2$ edges, and $3\cdot2<10$.

\emph{$\ell\in\{4,\dots,7,9,10,11\}$:} Proposition~\ref{prop:woodall} applies with $n=\R{3}{\ell}$
from the table; then $\ell\le n\le2\ell-3$, and~\eqref{eq:basic-wood} reads
\[
\begin{array}{c|ccccccc}
(\ell,n) & (4,5) & (5,6) & (6,7) & (7,8) & (9,11) & (10,12) & (11,13)\\
\hline
n(n-1) & 20 & 30 & 42 & 56 & 110 & 132 & 156\\
3\bigl[\binom{\ell-1}2+\binom{n-\ell+2}2\bigr] & 18 & 27 & 39 & 54 & 102 & 126 & 153
\end{array}
\]

\emph{$\ell=8$:} suppose a $3$-colouring of $K_9$ has no colour-avoiding $C_8$. By
Corollary~\ref{cor:mono}, $\circum(G_j)\le7$, so $e(G_j)\le24$ by Lemma~\ref{lem:basic-nine} (a graph with $e\ge25$ would be $K_7$ and $K_3$ sharing a vertex, which has $24$ edges), and
by~\eqref{eq:basic-sum} $e(G_j)=24$ for all $j$. Every edge at a vertex $v$ lies in exactly two
$G_j$, so $\sum_j\deg_{G_j}(v)=16$; as no three elements of $\{2,6\}$ sum to $16$,
Lemma~\ref{lem:basic-nine} shows $\deg_{G_j}(v)=8$ for some $j$. This gives $9$ pairs $(v,j)$ with
$\deg_{G_j}(v)=8$, whereas by Lemma~\ref{lem:basic-nine} there are only $3$.
\end{proof}

\section{Cycles of length at least 11}\label{sec:t3all}

Throughout this section $\ell\ge11$,
\[
  d:=\fl{\frac{\ell+1}{6}},\qquad n:=\ell+d,\qquad \ell=6d+\delta\ \text{ with }\ \delta\in\{-1,0,1,2,3,4\},
\]
so $d\ge2$, $n=7d+\delta\ge13$ and $6d-1\le\ell\le6d+4$.

\begin{theorem}\label{thm:t3all}
For every integer $\ell\ge11$, $\R{3}{\ell}=\ell+\fl{(\ell+1)/6}$.
\end{theorem}

The lower bound is Proposition~\ref{prop:lowA}. A $3$-colouring of $K_n$ is \emph{bad} if
$\circum(G_j)\le\ell-1$ for every $j$; by Corollary~\ref{cor:mono} it has no colour-avoiding $C_\ell$
if and only if it is bad. We show that no bad colouring of $K_n$ exists. The idea is to close two of
the three graphs $G_j$ (Bondy--Chv\'atal), which forces a clique-plus-attachment structure, to bound
the edges between the attached sets and the cliques by an insertion argument, and to compare with an
exact count of pairs.

\emph{Conventions.} $V$ is the vertex set of $K_n$, all graphs have vertex set $V$, and
$\{\alpha,\beta,\gamma\}=[3]$. A \emph{pair} is a $2$-subset of $V$. A pair of colour $\alpha$ lies
in $G_\beta$ and $G_\gamma$ but not in $G_\alpha$, so $\sum_j\deg_{G_j}(v)=2(n-1)$ for every $v$. Let
$c_j(v):=n-1-\deg_{G_j}(v)$, the number of pairs of colour $j$ at $v$, and for a graph $H$ let
$\mathrm{Mid}(H):=\{v:\deg_H(v)\ge n/2\}$.

\subsection{Closure and the normal form}

Let $\mathrm{cl}(H)$ be the closure of Theorem~\ref{thm:closure} with $N=n$. A graph $H$ is
\emph{closed} if $\deg_H(u)+\deg_H(v)\le n-1$ for all non-adjacent $u\ne v$; $\mathrm{cl}(H)$ is
closed.

\begin{lemma}\label{t3a:lem-cl}
$\circum(\mathrm{cl}(H))=\circum(H)$ for every graph $H$ on $V$.
\end{lemma}

\begin{proof}
It suffices to show $\circum(H+uv)\le\circum(H)$ when $u,v$ are non-adjacent with
$\deg_H(u)+\deg_H(v)\ge n$. Let $Z$ be a longest cycle of $H+uv$; if $uv\notin E(Z)$, then $Z\subseteq H$.
Otherwise $Z-uv=p_1p_2\cdots p_p$ with $p_1=u$, $p_p=v$, $p=|Z|$. Let
$S:=\{i\in[1,p-1]:up_{i+1}\in E(H)\}$ and $T:=\{i\in[1,p-1]:vp_i\in E(H)\}$. If $i\in S\cap T$, then
$p_1\cdots p_i\,p_pp_{p-1}\cdots p_{i+1}\,p_1$ is a cycle of length $p$ in $H$. Otherwise
$|S|+|T|\le p-1$, so $u$ and $v$ have at least $n-p+1$ neighbours off the path, counted with
multiplicity, among $n-p$ vertices; a common neighbour $w$ gives the cycle $p_1\cdots p_p\,w\,p_1$ of
length $p+1$ in $H$.
\end{proof}

\begin{lemma}\label{t3a:lem-closed}
Let $H$ be closed with $\circum(H)\le\ell-1$, and $Z:=\mathrm{Mid}(H)$.
\begin{enumerate}
\item[(i)] $Z$ is a clique of $H$, and $|Z|\le\ell-1$.
\item[(ii)] Every $x\notin Z$ has $\deg_H(x)\le n-|Z|$.
\item[(iii)] If $n$ is odd and $x\notin Z$ has $\deg_H(x)=(n-1)/2$, then $Z\subseteq N_H(x)$; in
  particular $|Z|\le(n-1)/2$.
\end{enumerate}
\end{lemma}

\begin{proof}
(i) Two vertices of $Z$ have degree sum at least $n$, so they are adjacent; a clique on $\ell$
vertices contains $C_\ell$. (ii) If $|Z|<n/2$, then $\deg_H(x)<n/2<n-|Z|$. Otherwise $\deg_H(x)<n/2\le|Z|$, so
$x$ has a non-neighbour $w\in Z$, $\deg_H(w)\ge|Z|-1$, and closedness gives $\deg_H(x)\le n-|Z|$. (iii) Every
$w\in Z$ has $\deg_H(w)\ge(n+1)/2$, so $xw\in E(H)$.
\end{proof}

\begin{lemma}\label{t3a:lem-nf}
Let $\phi$ be bad and $\alpha$ a colour. Let $\phi^\alpha$ be obtained from $\phi$ by giving colour
$\alpha$ to every pair of colour $j\in\{\beta,\gamma\}$ that is an edge of $\mathrm{cl}(G_j(\phi))$.
Then $G_j(\phi^\alpha)=\mathrm{cl}(G_j(\phi))$ for $j\in\{\beta,\gamma\}$,
$G_\alpha(\phi^\alpha)\subseteq G_\alpha(\phi)$, and $\phi^\alpha$ is bad. We call $\alpha$ the
\emph{sink}.
\end{lemma}

\begin{proof}
The recoloured pairs of colour $\beta$ are exactly the pairs of $E(\mathrm{cl}(G_\beta(\phi)))
\setminus E(G_\beta(\phi))$; such a pair enters $G_\beta$, leaves $G_\alpha$ and stays in $G_\gamma$.
The same holds with $\beta,\gamma$ exchanged. This gives the first two assertions, and badness
follows from Lemma~\ref{t3a:lem-cl} and $G_\alpha(\phi^\alpha)\subseteq G_\alpha(\phi)$.
\end{proof}

\begin{lemma}\label{t3a:lem-two}
Let $\phi$ be bad, $M_j:=\mathrm{Mid}(\mathrm{cl}(G_j(\phi)))$, and $t(v)$ the number of $j$ with
$v\in M_j$. Then $t(v)\ge2$ for every $v\in V$.
\end{lemma}

\begin{proof}
Write $G_j=G_j(\phi)$ and $\widehat G_j=\mathrm{cl}(G_j)$; by Lemma~\ref{t3a:lem-cl},
Lemma~\ref{t3a:lem-closed} applies to $\widehat G_j$, so $|M_j|\le\ell-1<n$.

\emph{Step 1.} Let $v\notin M_i\cup M_{i'}$ for two colours $i,i'$, and let $j$ be the third. For
$h\in\{i,i'\}$, $\deg_{G_h}(v)\le\deg_{\widehat G_h}(v)\le\fl{(n-1)/2}$, and $\deg_{G_j}(v)\le n-1$.
As the three degrees sum to $2(n-1)$, we get $2(n-1)\le2\fl{(n-1)/2}+n-1$, so $n$ is odd and
\begin{equation}\label{t3a:eq-tight}
  \deg_{G_j}(v)=n-1,\qquad \deg_{G_h}(v)=\deg_{\widehat G_h}(v)=\tfrac{n-1}2\quad(h\in\{i,i'\}).
\end{equation}
Hence $t(v)=1$ ($v\in M_j$). Moreover (a) $|M_i|,|M_{i'}|\le(n-1)/2$ by
Lemma~\ref{t3a:lem-closed}(iii), and (b) $M_i\cap M_{i'}=\emptyset$: for $w\in M_i\cap M_{i'}$,
Lemma~\ref{t3a:lem-closed}(iii) puts $vw$ in $\widehat G_i$ and $\widehat G_{i'}$, hence, as
$N_{G_h}(v)=N_{\widehat G_h}(v)$ by~\eqref{t3a:eq-tight}, in $G_i$ and $G_{i'}$; then $vw$ has colour
$j$, contradicting $\deg_{G_j}(v)=n-1$.

\emph{Step 2.} Suppose $t(v)\le1$ for some $v$; by Step~1, $t(v)=1$, say $v\in M_1$, and $M_2\cap
M_3=\emptyset$, $|M_2|,|M_3|\le(n-1)/2$. Take $u\notin M_1$; as $M_2\cap M_3=\emptyset$, Step~1 applies to $u$, so $t(u)=1$, say $u\in M_2$, and
$M_1\cap M_3=\emptyset$, $|M_1|\le(n-1)/2$. Then $|M_1\cup M_2|\le n-1$, so some $w$ lies in $M_3$
only, and Step~1 at $w$ gives $M_1\cap M_2=\emptyset$. So $M_1,M_2,M_3$ partition $V$, and $t\equiv1$.

\emph{Step 3.} By~\eqref{t3a:eq-tight} at every vertex, the vertices of $M_j$ have degree $n-1$ in
$G_j$ and the others degree $(n-1)/2$. As $n\ge4$, some $|M_j|\ge2$; then $G_j$ is $2$-connected (for
every $u$, a vertex of $M_j\setminus\{u\}$ is adjacent to all of $G_j-u$) with minimum degree
$(n-1)/2$, and Theorem~\ref{thm:dirac} gives $\circum(G_j)\ge n-1\ge\ell$, contradicting badness.
\end{proof}

\begin{corollary}\label{t3a:cor-sink}
If $\phi$ is bad, there is a colour $\alpha$ such that, with $G_j:=G_j(\phi^\alpha)$,
$\mathrm{Mid}(G_\beta)\cup\mathrm{Mid}(G_\gamma)=V$ and $|\mathrm{Mid}(G_\beta)\cap
\mathrm{Mid}(G_\gamma)|\ge n/3$.
\end{corollary}

\begin{proof}
By Lemma~\ref{t3a:lem-nf}, $\mathrm{Mid}(G_j(\phi^\alpha))=M_j$ for $j\ne\alpha$. By
Lemma~\ref{t3a:lem-two} every vertex lies in $M_\beta\cup M_\gamma$ for every choice of $\alpha$, and
$\sum_{\alpha}|M_\beta\cap M_\gamma|=\sum_v\binom{t(v)}2\ge n$ (for each $\alpha$, $\beta,\gamma$ are the other two colours); take $\alpha$ with the largest term.
\end{proof}

\subsection{Structure of the normal form}

\begin{lemma}[]\label{t3a:lem-close}
Let $\phi$ be a $3$-colouring of $K_n$, $a$ a colour and $j,j'$ the other two. Let $X,S\subseteq V$ be
disjoint, $m:=|X|$, $s:=|S|\ge2$, $\lambda:=\ell+1-2s$, and suppose that each $z\in S$ lies in at most
one pair $zx$ with $x\in X$ and $\phi(zx)\ne a$, and that
\[
  \lambda\ge4,\qquad m\ge\lambda+\fl{\lambda/2}-1,\qquad m\ge\lambda+s+1 .
\]
Then $G_j$ or $G_{j'}$ contains $C_\ell$.
\end{lemma}

\begin{proof}
For $z\in S$ and $\eta\in\{j,j'\}$ let $M_\eta(z):=\{x\in X:\phi(zx)=\eta\}$; so for $x\in X$,
$zx\in E(G_\eta)$ if and only if $x\notin M_\eta(z)$, and $|M_j(z)|+|M_{j'}(z)|\le1$. All paths below
have their vertices in $X$; a \emph{$\lambda$-path} is a path on $\lambda$ vertices.

\emph{Assembly.} Let $x_1\cdots x_\lambda$ be a path of $G_\eta$ and $z_1,\dots,z_s$ an enumeration of
$S$ with $x_\lambda\notin M_\eta(z_1)$ and $x_1\notin M_\eta(z_s)$. Choose distinct
$y_1,\dots,y_{s-1}\in X\setminus\{x_1,\dots,x_\lambda\}$ successively with $y_t\notin M_\eta(z_t)\cup
M_\eta(z_{t+1})$; at most $|M_\eta(z_t)|+|M_\eta(z_{t+1})|+(t-1)\le s<m-\lambda$ of the $m-\lambda$
candidates are excluded. Then $x_1\cdots x_\lambda\,z_1y_1z_2y_2\cdots z_{s-1}y_{s-1}z_s\,x_1$ is a
cycle of $G_\eta$ of length $\lambda+2s-1=\ell$.

\emph{Paths.} Colour a pair inside $X$ red if its colour is not $j$ and blue otherwise; red pairs are
edges of $G_j$ and blue pairs edges of $G_{j'}$. By Theorem~\ref{thm:gg} with $a=b=\lambda$, every set of at least
$\lambda+\fl{\lambda/2}-1$ vertices of $X$ contains a $\lambda$-path of $G_j$ or of $G_{j'}$.

Fix distinct $w,w'\in S$. Call $(x,\eta)\in X\times\{j,j'\}$ \emph{dead} if $x\in M_\eta(w)\cap
M_\eta(w')$; as $|M_j(w)|+|M_{j'}(w)|\le1$, there is at most one dead pair. A path $x_1\cdots
x_\lambda$ of $G_\eta$ is \emph{good} if neither $(x_1,\eta)$ nor $(x_\lambda,\eta)$ is dead.

\emph{A good $\lambda$-path of $G_\eta$ gives $C_\ell\subseteq G_\eta$.} Take $u_1,u_\lambda\in
\{w,w'\}$ with $x_1\notin M_\eta(u_1)$, $x_\lambda\notin M_\eta(u_\lambda)$. If $u_1\ne u_\lambda$,
put $z_s:=u_1$, $z_1:=u_\lambda$. If $u_1=u_\lambda$, let $u'$ be the other vertex of $\{w,w'\}$; as
$|M_\eta(u')|\le1$, either $x_1\notin M_\eta(u')$ and we put $(z_1,z_s):=(u_1,u')$, or
$x_\lambda\notin M_\eta(u')$ and we put $(z_1,z_s):=(u',u_1)$. Enumerating the rest of $S$ as
$z_2,\dots,z_{s-1}$, the assembly applies.

It remains to find a good $\lambda$-path. If there is no dead pair, (Paths) gives one. Otherwise let
$(x^*,\eta^*)$ be the dead pair; the hypotheses are symmetric in $j,j'$, so we may assume $\eta^*=j$
(and define red/blue afterwards). Then all paths of $G_{j'}$, and all paths of $G_j$ not ending at
$x^*$, are good. Suppose there is no good $\lambda$-path. By (Paths) there is a $\lambda$-path, so it
is a path $x_1\cdots x_{\lambda-1}x^*$ of $G_j$. If $x^*$ had a $G_j$-neighbour $x_p$ with
$p\le\lambda-2$, then $x_1\cdots x_p\,x^*x_{\lambda-1}\cdots x_{p+1}$ would be good; if it had a
$G_j$-neighbour $w_0\in X$ off the path, then $x_2\cdots x_{\lambda-1}x^*w_0$ would be good. So
$x^*w_0\in E(G_{j'})$ for all $w_0\in X\setminus\{x^*,x_{\lambda-1}\}$. Hence $X\setminus\{x^*\}$
contains no path $w_1\cdots w_{\lambda-2}$ of $G_{j'}$: one end, say $w_1$, differs from
$x_{\lambda-1}$, and for any $w''\in X\setminus\{x^*,x_{\lambda-1},w_1,\dots,w_{\lambda-2}\}$ (which
exists as $m\ge\lambda+1$) the path $w''x^*w_1\cdots w_{\lambda-2}$ would be a good $\lambda$-path of
$G_{j'}$. Also $X\setminus\{x^*\}$ contains no $\lambda$-path of $G_j$. As $m-1\ge\lambda+
\fl{(\lambda-2)/2}-1$, this contradicts Theorem~\ref{thm:gg} with orders $\lambda\ge\lambda-2\ge2$.
\end{proof}

\begin{definition}[assumptions]\label{t3a:def-S}
Until Subsection~\ref{t3a:sub-end}, $\phi$ is bad, $\alpha$ is as in Corollary~\ref{t3a:cor-sink},
$G_\alpha,G_\beta,G_\gamma$ are the avoiding graphs of $\phi^\alpha$, and colours refer to
$\phi^\alpha$; so $\phi^\alpha$ is bad and $G_\beta,G_\gamma$ are closed. Put
$\mathrm{Mid}_j:=\mathrm{Mid}(G_j)$, $m_j:=|\mathrm{Mid}_j|$, $P_j:=V\setminus\mathrm{Mid}_j$,
$p_j:=|P_j|$ ($j\in\{\beta,\gamma\}$), $K:=\mathrm{Mid}_\beta\cap\mathrm{Mid}_\gamma$, $\kappa:=|K|$,
$Q:=V\setminus K$, $q:=|Q|$ and $\sigma:=\ell-1-q$. As $\mathrm{Mid}_\beta\cup\mathrm{Mid}_\gamma=V$,
the sets $K,P_\beta,P_\gamma$ partition $V$, $\mathrm{Mid}_\beta=K\cup P_\gamma$,
$\mathrm{Mid}_\gamma=K\cup P_\beta$, and
\[
  q=p_\beta+p_\gamma,\qquad m_\beta=n-p_\beta=\kappa+p_\gamma,\qquad m_\gamma=n-p_\gamma=\kappa+p_\beta,
  \qquad\sigma=\kappa-d-1 .
\]
For $j\in\{\beta,\gamma\}$, $B_j$ is the set of edges of $G_j$ between $P_j$ and $\mathrm{Mid}_j$,
and $h_j:=\min\{p_j,\fl{n/2}-1\}$. For $x\in K$, $b(x)$ is the number of pairs $xy$ of colour $\beta$
with $y\in P_\beta$, and $g(x)$ the number of pairs $xy$ of colour $\gamma$ with $y\in P_\gamma$.
\end{definition}

\begin{lemma}\label{t3a:lem-struct}
\begin{enumerate}
\item[(a)] Pairs inside $\mathrm{Mid}_\gamma$ have colour $\alpha$ or $\beta$, pairs inside
  $\mathrm{Mid}_\beta$ colour $\alpha$ or $\gamma$, and pairs inside $K$ colour $\alpha$.
\item[(b)] $p_\beta,p_\gamma\ge d+2$.
\item[(c)] $3\kappa\ge n$; in particular $\kappa\ge d+1\ge3$ and $\sigma\ge0$.
\item[(d)] For $x\in K$, $N_{G_\alpha}(x)\subseteq Q$, $\deg_{G_\alpha}(x)=b(x)+g(x)$, $b(x)\le
  h_\beta$ and $g(x)\le h_\gamma$.
\end{enumerate}
\end{lemma}

\begin{proof}
(a) $\mathrm{Mid}_j$ is a clique of $G_j$ by Lemma~\ref{t3a:lem-closed}(i).

(b) Let $j\in\{\beta,\gamma\}$. By Lemma~\ref{t3a:lem-closed}(i), $m_j\le\ell-1$; suppose
$m_j=\ell-1$ and let $S:=P_j$, so $s:=|S|=d+1\ge3$. A vertex $z\in S$ with two $G_j$-neighbours
$u\ne v$ in $\mathrm{Mid}_j$ would close a Hamilton $u$--$v$ path of the clique $\mathrm{Mid}_j$ to a
$C_\ell$ in $G_j$; so $z$ lies in at most one pair $zx$, $x\in\mathrm{Mid}_j$, of colour $\ne j$.
Apply Lemma~\ref{t3a:lem-close} to $\phi^\alpha$ with $a:=j$, $X:=\mathrm{Mid}_j$, $m=\ell-1$, and
$\lambda=\ell-1-2d$. Its hypotheses hold: $\lambda\ge(2\ell-4)/3\ge6$ as $2d\le(\ell+1)/3$;
$m\ge\lambda+\fl{\lambda/2}-1$ is equivalent to $\fl{\lambda/2}\le2d+1$, true as $\lambda\le4d+3$;
and $m\ge\lambda+s+1$ is equivalent to $d\ge2$. So $\phi^\alpha$ is not bad, a contradiction. Hence
$m_j\le\ell-2$, that is, $p_j=n-m_j\ge d+2$.

(c) $3\kappa\ge n$ is Corollary~\ref{t3a:cor-sink}, and $n/3\ge d+1\ge3$ as $\ell\ge2d+3$.

(d) By (a), the edges of $G_\alpha$ at $x$, that is, the pairs of colour $\beta$ or $\gamma$, go to
$P_\beta$ (colour $\beta$) or $P_\gamma$ (colour $\gamma$). Clearly $b(x)\le p_\beta$, and as
$x\in\mathrm{Mid}_\beta$, $b(x)\le c_\beta(x)\le n-1-\cl{n/2}=\fl{n/2}-1$; likewise for $g(x)$.
\end{proof}

\subsection{Insertion}

Let $\mathcal N=(N(x))_{x\in X}$ be a family of subsets of a finite set $Z$. A \emph{rainbow linear
forest} of $\mathcal N$ is a set $Y\subseteq X$ with $2$-sets $\pi(x)\subseteq N(x)$ ($x\in Y$),
pairwise distinct, such that the graph $F$ on $Z$ with edge set $\{\pi(x):x\in Y\}$ is a linear forest
(a disjoint union of paths). Its size is $|Y|$; $\mathrm{lf}(\mathcal N)$ is the largest size. A
vertex is \emph{free} if its degree in $F$ is at most $1$. A subset of a rainbow linear forest is
one.

\begin{lemma}\label{t3a:lem-greedy}
If $(Y,\pi)$ is a rainbow linear forest of size $k$ and $x\in X\setminus Y$ has $|N(x)|\ge k+2$, then
$\pi(x)$ can be chosen so that $(Y\cup\{x\},\pi)$ is a rainbow linear forest.
\end{lemma}

\begin{proof}
A path with $i\ge1$ edges has $i-1$ vertices of degree $2$, so $F$ has at most $\max\{k-1,0\}$
non-free vertices, and $N(x)$ contains at least two free vertices if $k=0$ and three if $k\ge1$. A
component contains at most two free vertices (and one if $k=0$), so $N(x)$ contains free $u,w$ in
different components, and $\pi(x):=\{u,w\}$ works.
\end{proof}

\begin{theorem}\label{t3a:thm-inse}
Let $r:=|X|$, $D:=\max_x|N(x)|$ ($D:=0$ if $X=\emptyset$) and $e:=\sum_x|N(x)|$. Then
$e\le(\mathrm{lf}(\mathcal N)+1)\max\{r,D\}$.
\end{theorem}

\begin{proof}
Induction on $r$, trivial for $r=0$. Let $f:=\mathrm{lf}(\mathcal N)$; call $x$ \emph{rich} if
$|N(x)|\ge f+2$, let $q'$ be the number of rich indices and $\mathcal N_p$ the subfamily of the other
(poor) indices. If $q'=0$, then $e\le(f+1)r$. Otherwise, starting from a largest rainbow linear forest
of $\mathcal N_p$, of size $f_p$, and adding the rich indices one by one by
Lemma~\ref{t3a:lem-greedy} (each intermediate forest has size at most $f$) gives $f_p+q'\le f$. By
induction, as $\mathcal N_p$ has $r-q'$ members of size at most $f+1$, and as $D\ge f+2$,
\[
  e\le(f_p+1)\max\{r-q',f+1\}+q'D\le(f-q'+1)\max\{r,D\}+q'\max\{r,D\}=(f+1)\max\{r,D\}. \qedhere
\]
\end{proof}

\begin{lemma}\label{t3a:lem-insert}
Let $\circum(H)\le\ell-1$, let $Z$ be a clique of $H$ with $|Z|\ge3$, $X\subseteq V(H)\setminus Z$
and $N(x):=N_H(x)\cap Z$. Then $\mathrm{lf}((N(x))_{x\in X})\le\ell-1-|Z|$.
\end{lemma}

\begin{proof}
Let $(Y,\pi)$ be a rainbow linear forest of size $k$. Listing the components of $F$ one after the
other orders $Z$ as $z_1,\dots,z_{|Z|}$ so that every edge of $F$ is some $z_iz_{i+1}$. Replacing each
edge $\pi(x)=z_iz_{i+1}$ of the cycle $z_1z_2\cdots z_{|Z|}z_1$ of $H$ by the path $z_ixz_{i+1}$
gives a cycle of length $|Z|+k$ in $H$ (the $x\in Y$ are distinct vertices outside $Z$ and the
$\pi(x)$ are distinct edges), so $|Z|+k\le\ell-1$.
\end{proof}

\begin{corollary}\label{t3a:cor-B}
$|B_j|\le(p_j-d)\,p_j$ for $j\in\{\beta,\gamma\}$.
\end{corollary}

\begin{proof}
Let $N(x):=N_{G_j}(x)\cap\mathrm{Mid}_j$ for $x\in P_j$; then $\sum_x|N(x)|=|B_j|$, $r=p_j$, and
$|N(x)|\le n-m_j=p_j$ by Lemma~\ref{t3a:lem-closed}(ii). Lemma~\ref{t3a:lem-insert} with $H=G_j$,
$Z=\mathrm{Mid}_j$ (a clique with $m_j\ge\kappa\ge3$) gives $\mathrm{lf}\le\ell-1-m_j$, and
Theorem~\ref{t3a:thm-inse} gives $|B_j|\le(\ell-m_j)p_j=(p_j-d)p_j$.
\end{proof}

\subsection{The sink colour}

\begin{lemma}\label{t3a:lem-qh}
Let $\widehat G:=\mathrm{cl}(G_\alpha)$. Then $\circum(\widehat G)\le\ell-1$ and $Q$ is a clique of
$\widehat G$.
\end{lemma}

\begin{proof}
The first assertion is Lemma~\ref{t3a:lem-cl}. For $y\in P_\beta$ we have
$\deg_{G_\beta}(y)\le\cl{n/2}-1$, so $c_\beta(y)\ge\fl{n/2}$ and
\begin{equation}\label{t3a:eq-qh}
  \deg_{G_\alpha}(y)=c_\beta(y)+c_\gamma(y)\ge\fl{n/2};
\end{equation}
symmetrically for $z\in P_\gamma$, with $c_\gamma(z)\ge\fl{n/2}$. Suppose $y,y'\in Q$ are not
adjacent in $\widehat G$. As $\widehat G$ is closed and contains $G_\alpha$, \eqref{t3a:eq-qh} forces
$n$ odd and $\deg_{\widehat G}(y)=\deg_{G_\alpha}(y)=(n-1)/2$; by symmetry between $\beta$ and
$\gamma$ let $y\in P_\beta$. Then $c_\beta(y)=(n-1)/2$, $c_\gamma(y)=0$, $\deg_{G_\beta}(y)=(n-1)/2$,
and Lemma~\ref{t3a:lem-closed}(iii) gives $m_\beta\le(n-1)/2$ and $y$ adjacent in $G_\beta$ to all of
$\mathrm{Mid}_\beta$. Now $P_\gamma\ne\emptyset$, since otherwise $\mathrm{Mid}_\gamma=V$ would be a
clique on $n\ge\ell$ vertices. For $z\in P_\gamma\subseteq\mathrm{Mid}_\beta$, the pair $yz$ lies in
$G_\beta$ and does not have colour $\gamma$ (as $c_\gamma(y)=0$), so it has colour $\alpha$ and is not
an edge of $\widehat G$ (as $N_{\widehat G}(y)=N_{G_\alpha}(y)$). Closedness gives
$\deg_{\widehat G}(z)\le(n-1)/2$, so by~\eqref{t3a:eq-qh} for $z$, $c_\gamma(z)=(n-1)/2$,
$c_\beta(z)=0$ and $\deg_{G_\gamma}(z)=(n-1)/2$, and Lemma~\ref{t3a:lem-closed}(iii) gives
$m_\gamma\le(n-1)/2$. Then $m_\beta+m_\gamma\le n-1$, whereas $m_\beta+m_\gamma=n+\kappa$.
\end{proof}

\begin{proposition}\label{t3a:prop-scl}
Let $\widehat G:=\mathrm{cl}(G_\alpha)$ and $2q\ge\ell+1$. Let $K^+:=\{x\in K:Q\subseteq
N_{\widehat G}(x)\}$, $k^+:=|K^+|$ and $K^-:=K\setminus K^+$. Then
\begin{enumerate}
\item[(a)] every $x\in K^-$ has $\deg_{\widehat G}(x)\le\kappa-k^+$;
\item[(b)] $k^+\le\sigma$, and $\mathrm{lf}(\mathcal N^-)\le\sigma-k^+$ for $\mathcal
  N^-:=(N_{\widehat G}(x)\cap Q)_{x\in K^-}$;
\item[(c)] $e_{G_\alpha}(K^-,Q)\le(\sigma+1-k^+)(\kappa-k^+)$.
\end{enumerate}
\end{proposition}

\begin{proof}
(a) A non-neighbour $y\in Q$ of $x$ in $\widehat G$ is adjacent to the other $q-1$ vertices of $Q$
(Lemma~\ref{t3a:lem-qh}) and to $K^+$, so closedness gives $\deg_{\widehat G}(x)\le
n-1-(q-1+k^+)=\kappa-k^+$.

(b) Let $f^-:=\mathrm{lf}(\mathcal N^-)$ and suppose $k^++f^-\ge s:=\sigma+1$. Take a rainbow linear
forest of $\mathcal N^-$ of size $\min\{f^-,s\}$ and, while its size $k$ is less than $s$, add an
unused index $x\in K^+$ by Lemma~\ref{t3a:lem-greedy}, with $N(x)=Q$ and $|Q|=q\ge k+2$, as $k\le
s-1=\ell-1-q\le q-2$ by $2q\ge\ell+1$. At most $k^+$ indices are needed, giving a rainbow linear forest
of $(N_{\widehat G}(x)\cap Q)_{x\in K}$ of size $\sigma+1=\ell-q$. This contradicts
Lemma~\ref{t3a:lem-insert} with $H=\widehat G$ and the clique $Z=Q$ ($q\ge3$ by
Lemma~\ref{t3a:lem-struct}(b)).

(c) By (a), $\mathcal N^-$ has $\kappa-k^+$ members of size at most $\kappa-k^+$, so
Theorem~\ref{t3a:thm-inse} and (b) give $e_{\widehat G}(K^-,Q)\le(\sigma+1-k^+)(\kappa-k^+)$; and
$G_\alpha\subseteq\widehat G$.
\end{proof}

\subsection{The count}

\begin{lemma}\label{t3a:lem-count}
With $\mathrm{Bb}:=\sum_{x\in K}b(x)$ and $\mathrm{Cc}:=\sum_{x\in K}g(x)$,
\[
  \mathrm{Bb}+|B_\beta|\le m_\beta p_\beta,\qquad \mathrm{Cc}+|B_\gamma|\le m_\gamma p_\gamma,\qquad
  \mathrm{Bb}+\mathrm{Cc}+|B_\beta|+|B_\gamma|\ge\kappa q+p_\beta p_\gamma .
\]
\end{lemma}

\begin{proof}
Of the $m_\beta p_\beta$ pairs between $P_\beta$ and $\mathrm{Mid}_\beta=K\cup P_\gamma$, those not of
colour $\beta$ form $B_\beta$, and those of colour $\beta$ are the $\mathrm{Bb}$ pairs to $K$ and some
number $N_\beta$ of pairs to $P_\gamma$. So $|B_\beta|=m_\beta p_\beta-\mathrm{Bb}-N_\beta$, and
likewise $|B_\gamma|=m_\gamma p_\gamma-\mathrm{Cc}-N_\gamma$, where $N_\beta+N_\gamma\le p_\beta
p_\gamma$. Finally $m_\beta p_\beta+m_\gamma p_\gamma-p_\beta p_\gamma=\kappa q+p_\beta p_\gamma$.
\end{proof}

For integers $p_\beta,p_\gamma,k$ put $q:=p_\beta+p_\gamma$, $\kappa:=n-q$, $m_j:=n-p_j$,
$h_j:=\min\{p_j,\fl{n/2}-1\}$, and
\begin{gather*}
  \mathrm{w}(p_\beta,p_\gamma):=p_\beta^2-p_\beta p_\gamma+p_\gamma^2-dq,\\
  \mathrm{mg}(p_\beta,p_\gamma,k):=\sum_{j\in\{\beta,\gamma\}}\min\bigl\{(p_j-d)p_j+k\,h_j,\;
  m_jp_j\bigr\}+(\kappa-d-k)(\kappa-k)-\kappa q-p_\beta p_\gamma .
\end{gather*}

\begin{theorem}\label{t3a:thm-hv}
\begin{enumerate}
\item[(i)] $\mathrm{w}(p_\beta,p_\gamma)\ge0$, and $\mathrm{w}(p_\beta,p_\gamma)\ge1$ if $2q=\ell$.
\item[(ii)] If $2q\ge\ell+1$, then $k:=k^+$ satisfies $0\le k\le\sigma$ and
  $\mathrm{mg}(p_\beta,p_\gamma,k)\ge0$.
\end{enumerate}
\end{theorem}

\begin{proof}
By Lemma~\ref{t3a:lem-struct}(d), $\mathrm{Bb}+\mathrm{Cc}=e_{G_\alpha}(K,Q)$.

(i) By Lemma~\ref{t3a:lem-count} and Corollary~\ref{t3a:cor-B}, $\kappa q+p_\beta p_\gamma\le
e_{G_\alpha}(K,Q)+(p_\beta-d)p_\beta+(p_\gamma-d)p_\gamma$, and $e_{G_\alpha}(K,Q)\le\kappa q$. If
$2q=\ell$ and $e_{G_\alpha}(K,Q)=\kappa q$, then, as $\kappa=q+d>q\ge2$, the complete bipartite graph
between $K$ and $Q$ in $G_\alpha$ contains $C_{2q}=C_\ell$, contradicting badness. This gives (i).

(ii) $0\le k\le\sigma$ by Proposition~\ref{t3a:prop-scl}(b). Split $\mathrm{Bb}=B'+B''$ and
$\mathrm{Cc}=C'+C''$ into the sums over $K^+$ and over $K^-$. By Lemma~\ref{t3a:lem-struct}(d),
$B'\le kh_\beta$ and $C'\le kh_\gamma$, and by Proposition~\ref{t3a:prop-scl}(c), with
$\sigma+1=\kappa-d$, $B''+C''=e_{G_\alpha}(K^-,Q)\le(\kappa-d-k)(\kappa-k)$. By
Lemma~\ref{t3a:lem-count}, Corollary~\ref{t3a:cor-B} and $\min\{u+w,v\}\le\min\{u,v\}+w$ ($w\ge0$),
\[
  \kappa q+p_\beta p_\gamma\le\sum_{j\in\{\beta,\gamma\}}\min\bigl\{(p_j-d)p_j+kh_j,\,m_jp_j\bigr\}
  +B''+C'',
\]
which gives $\mathrm{mg}(p_\beta,p_\gamma,k)\ge0$.
\end{proof}

\subsection{The arithmetic}

\begin{proposition}\label{t3a:prop-arith}
Let $p_\beta,p_\gamma\ge d+2$ be integers, $q:=p_\beta+p_\gamma$ and $\kappa:=n-q$.
\begin{enumerate}
\item[(i)] If $2q\le\ell$, then $\mathrm{w}(p_\beta,p_\gamma)\le0$, with equality only if $2q=\ell$.
\item[(ii)] If $2q\ge\ell+1$ and $3\kappa\ge n$, then $\mathrm{mg}(p_\beta,p_\gamma,k)<0$ for every
  integer $k\in[0,\kappa-d-1]$.
\end{enumerate}
\end{proposition}

\begin{proof}
(i) With $x:=p_\beta-d-2\ge0$ and $y:=p_\gamma-d-2\ge0$,
\[
  \mathrm{w}(p_\beta,p_\gamma)=4-d^2+2(x+y)+x^2-xy+y^2\le4-d^2+2(x+y)+(x+y)^2 .
\]
As $x+y=q-2d-4\le\fl{\ell/2}-2d-4\le d-2$ and the right-hand side increases in $x+y\ge0$,
$\mathrm{w}\le4-d^2+2(d-2)+(d-2)^2=4-2d\le0$. Equality forces $d=2$ and $x=y=0$, so $q=8$ and
$16=2q\le\ell\le6d+4=16$.

(ii) \emph{Relaxation.} Put $a:=p_\beta$, $b:=p_\gamma$, $s:=\kappa-k\in[d+1,\kappa]$; as
$\mathrm{mg}$ is symmetric in $p_\beta,p_\gamma$, let $a\le b$. As $k\ge0$ and $h_j\le p_j$,
$\mathrm{mg}\le\overline{\mathrm{mg}}:=a\min\{a-d+k,n-a\}+b\min\{b-d+k,n-b\}+s(s-d)-\kappa q-ab$. We
show $\overline{\mathrm{mg}}<0$, treating all variables as reals.

\emph{Step 1: closed forms.} As $n+d-k=q+d+s$, we have $a-d+k\le n-a$ if and only if $a-b\le d+s$
(true), and $b-d+k\le n-b$ if and only if $b-a\le d+s$. If $b-a\le d+s$, then
\[
  \overline{\mathrm{mg}}=V(a,b,s):=a^2-ab+b^2-(s+d)q+s^2-ds=q^2-3ab-(s+d)q+s^2-ds .
\]
If $b-a\ge d+s$, then, as $n-b=\kappa+a$,
$\overline{\mathrm{mg}}=W(a,s):=a^2-as+s^2-d(a+s)$. Moreover $V(a,a+d+s,s)=W(a,s)$.

\emph{Step 2: $W<0$ on $\Omega:=\{(a,s):a\ge d+2,\ s\ge d+1,\ 2a+d+s\le2n/3\}$.} With $x:=a-d-2\ge0$,
$y:=s-d-1\ge0$,
\[
  W(a,s)=3-d^2+3x+x^2-xy+y^2,\qquad 2x+y\le R:=\tfrac{2n}{3}-4d-5=\tfrac{2d+2\delta-15}{3}.
\]
If $\Omega\ne\emptyset$, then $R\ge0$, so $2d\ge15-2\delta\ge7$, $d\ge4$, and $x\le R/2$, $y\le R$,
$R\le(2d-7)/3<2d/3$. Hence $W\le3-d^2+\frac32R+\frac54R^2<3+d-\frac49d^2\le3-\frac79d<0$, using
$d^2\ge4d$.

\emph{Step 3.} If $b-a\ge d+s$, then $2a+d+s\le q\le2n/3$ (as $3\kappa\ge n$), so
$(a,s)\in\Omega$ and $\overline{\mathrm{mg}}=W(a,s)<0$.

\emph{Step 4.} Let $b-a\le d+s$. For fixed $q,s$, $V$ decreases in $ab$; as $a'\le a\le q/2$ for
$a':=\max\{d+2,(q-d-s)/2\}$ (the second term as $q-2a=b-a\le d+s$), we get $ab\ge a'(q-a')$ (as $t\mapsto t(q-t)$ increases for $t\le q/2$) and
$V(a,b,s)\le V(a',q-a',s)$. If $a'=(q-d-s)/2$, then $(q-a')-a'=d+s$, so $V(a',q-a',s)=W(a',s)$ by Step~1 with $(a',s)\in\Omega$
(as $2a'+d+s=q\le2n/3$), which is negative. Otherwise $a'=d+2$, $s\ge q-3d-4$, and
\[
  \overline{\mathrm{mg}}\le V_0(q,s):=V(d+2,q-d-2,s)=q^2-3(d+2)(q-d-2)-(s+d)q+s^2-ds,
\]
where $q_1\le q\le2n/3$ with $q_1:=\fl{\ell/2}+1$ (as $2q\ge\ell+1$), and $\max\{d+1,q-3d-4\}\le s\le n-q$. As $V_0$ is
convex in $s$, it suffices to bound it at the two ends of this interval:
\begin{itemize}
\item $s=n-q$: $V_0(q,n-q)=\varphi(q):=3q^2-3(n+d+2)q+3(d+2)^2+n^2-dn$, convex in $q$, so
  $V_0\le\max\{\varphi(q_1),\varphi(2n/3)\}$;
\item $s=d+1$ and $q\le4d+5$: $V_0(q,d+1)=\psi(q):=q^2-(5d+7)q+3(d+2)^2+d+1$, convex, and $q_1\le
  q\le4d+5$ ($q_1\le3d+3$), so $V_0\le\max\{\psi(q_1),\psi(4d+5)\}$;
\item $s=q-3d-4$ and $q\ge4d+5$: $(q-d-2)-(d+2)=d+s$, so $V_0=W(d+2,s)$, and $(d+2,s)\in\Omega$ as $s\ge d+1$ and
  $2(d+2)+d+s=q\le2n/3$; so $V_0<0$.
\end{itemize}

\emph{Step 5: the values.} We have $q_1=3d+e$ with $e:=\fl{\delta/2}+1$, so $e=0,1,1,2,2,3$ for
$\delta=-1,0,1,2,3,4$ respectively, and with $n=7d+\delta$,
\begin{align*}
  \varphi(q_1)&=(4\delta-6e-6)\,d+\delta^2-3\delta e+3e^2-6e+12,\\
  \varphi(2n/3)&=-\tfrac53d^2+\bigl(\tfrac53\delta-16\bigr)d+\tfrac13(\delta-6)^2,\\
  \psi(q_1)&=-3d^2+(e-8)d+e^2-7e+13,\qquad \psi(4d+5)=3-d^2 .
\end{align*}
For $\delta=-1,0,1,2,3,4$, $\varphi(q_1)$ equals $13-10d$, $9-12d$, $7-8d$, $4-10d$, $3-6d$, $1-8d$,
negative for $d\ge2$. As $-1\le\delta\le4$, $\varphi(2n/3)\le-\frac53d^2-\frac{28}3d+\frac{49}3<0$. As
$0\le e\le3$, $\psi(q_1)\le-3d^2-5d+13<0$. Finally $3-d^2<0$.
\end{proof}

\subsection{Proof of Theorem~\ref{thm:t3all}}\label{t3a:sub-end}

Suppose $\phi$ is a bad $3$-colouring of $K_n$, choose $\alpha$ by Corollary~\ref{t3a:cor-sink}, and
use Definition~\ref{t3a:def-S}. By Lemma~\ref{t3a:lem-struct}(b),(c), $p_\beta,p_\gamma\ge d+2$ and
$3\kappa\ge n$. If $2q\le\ell$, Theorem~\ref{t3a:thm-hv}(i) and Proposition~\ref{t3a:prop-arith}(i)
give $\mathrm{w}(p_\beta,p_\gamma)=0$ and $2q=\ell$, and then Theorem~\ref{t3a:thm-hv}(i) gives
$\mathrm{w}(p_\beta,p_\gamma)\ge1$, a contradiction. If $2q\ge\ell+1$, the integer $k$ of
Theorem~\ref{t3a:thm-hv}(ii) satisfies $0\le k\le\sigma=\kappa-d-1$ and $\mathrm{mg}\ge0$,
contradicting Proposition~\ref{t3a:prop-arith}(ii). So no bad colouring exists, every $3$-colouring
of $K_n$ has a colour-avoiding $C_\ell$, and $\R{3}{\ell}\le n$. With Proposition~\ref{prop:lowA}
this proves Theorem~\ref{thm:t3all}.

\providecommand{\bysame}{\leavevmode\hbox to3em{\hrulefill}\thinspace}
\providecommand{\MR}{\relax\ifhmode\unskip\space\fi MR }
\providecommand{\MRhref}[2]{%
  \href{http://www.ams.org/mathscinet-getitem?mr=#1}{#2}
}
\providecommand{\href}[2]{#2}

\end{document}